\documentclass[10pt]{article}

\usepackage{amssymb,amsbsy,amsmath,amsfonts,amssymb,amscd}
\usepackage{latexsym,euscript,exscale,epic,eepic,epsfig}
\usepackage{latexsym,euscript,exscale,epic,eepic,epsfig}
\usepackage[english,francais]{babel}
\usepackage{amssymb}

\newtheorem{theorem}{Theorem}[section]
\newtheorem{lemma}[theorem]{Lemma}
\newtheorem{proposition}[theorem]{Proposition}
\newtheorem{corollary}[theorem]{Corollary}
\newtheorem{definition}[theorem]{Definition\rm}
\newtheorem{remark}{\it Remark\/}

\def\og{\leavevmode\raise.3ex\hbox{$\scriptscriptstyle\langle\!\langle$~}}
\def\fg{\leavevmode\raise.3ex\hbox{~$\!\scriptscriptstyle\,\rangle\!\rangle$}}
\begin{document}
\selectlanguage{english}
\vskip 0.5\baselineskip
\title{Multipliers spaces and pseudo-differential operators}
\author{By Sadek Gala\\ Universit\'e d'Evry
Val d'Essonne\\ D\'epartement de math\'ematiques\\Bd
F. Mitterrand. 91025 Evry Cedex. France\\Sadek.Gala@maths.univ-evry.fr }
\maketitle

\begin{abstract} Let $\sigma\left(  x,\xi\right)  $ be a sufficiently regular function defined
on $\mathbb{R}^{d}\times\mathbb{R}^{d}.$ The pseudo-differential operator with
symbol $\sigma$ is defined on the Schwartz class by the formula :%
\[
f\rightarrow\sigma f(x)=%
{\displaystyle\int\limits_{\mathbb{R}^{d}}}
\sigma\left(  x,\xi\right)  \widehat{f}(\xi)e^{2\pi ix\xi}d\xi,
\]
where $\widehat{f}(\xi)=%
{\displaystyle\int\limits_{\mathbb{R}^{d}}}
f(x)e^{-2\pi ix\xi}dx$ is the Fourier transform of $f.$\\
In this paper, we shall consider the regularity of the following type :

\begin{description}
\item[(a)] $\left\vert \partial_{\xi}^{\alpha}\sigma\left(  x,\xi\right)
\right\vert \leq A_{\alpha}\left(  1+\left\vert \xi\right\vert \right)
^{-\left\vert \alpha\right\vert },$

\item[(b)] $\left\vert \partial_{\xi}^{\alpha}\sigma\left(  x+y,\xi\right)
-\partial_{\xi}^{\alpha}\sigma\left(  x,\xi\right)  \right\vert \leq
A_{\alpha}\omega\left(  \left\vert y\right\vert \right)  \left(  1+\left\vert
\xi\right\vert \right)  ^{-\left\vert \alpha\right\vert },$

\end{description}

$\left(  \alpha\in\mathbb{N}^{d}\right)  $ where $\omega$ is suitable positive
function and we prove boundedness results for pseudo-differential operators on
multipliers spaces $X^{r}=\mathcal{M}\left(  H^{r}\rightarrow L^{2}\right)  $
whose symbol $\sigma\left(  x,\xi\right)  $ satisfies the regularity
condition on $x$.
\vskip 0.5\baselineskip
\noindent



\end{abstract}

\section{Introduction}

A pseudo-differential operator $\sigma$ with symbol $\sigma\left(
x,\xi\right)  $, defined initially on the Schwartz class of testing functions
$\mathcal{S}\left(  \mathbb{R}^{d}\right)  $, is given by%
\begin{equation}
f\rightarrow\sigma f(x)=%
{\displaystyle\int\limits_{\mathbb{R}^{d}}}
\sigma\left(  x,\xi\right)  \widehat{f}(\xi)e^{2\pi ix\xi}d\xi, \label{eq1}%
\end{equation}
with
\[
\widehat{f}(\xi)=%
{\displaystyle\int\limits_{\mathbb{R}^{d}}}
f(x)e^{-2\pi ix\xi}dx.
\]
being the Fourier transform of $f$.

We shall consider the standard symbol class, denoted by
$\boldsymbol{S}^{m}$, which is the most common and useful of the general
symbol classes. A function $\sigma$ belongs to $\boldsymbol{S}^{m}$ $\left(
\text{and is said to be of order }m\right)  $ if $\sigma\left(  x,\xi\right)
$ is a $C^{\infty}$ function of $\left(  x,\xi\right)  \in\mathbb{R}^{d}%
\times\mathbb{R}^{d}$ and satisfies the differential inequalities
\[
\left\vert \partial_{\xi}^{\alpha}\partial_{\xi}^{\beta}\sigma\left(
x,\xi\right)  \right\vert \leq A_{\alpha,\beta}\left(  1+\left\vert
\xi\right\vert \right)  ^{m-\left\vert \alpha\right\vert },
\]
for all multi-indices $\alpha$ and $\beta$.

Before we state our result, we need to make precise the definition of the
pseudo-differential operator (\ref{eq1}) and the class of symbols that is used.

We call the modulus of continuity every function $\omega:\left[
0,\mathbb{+\infty}\right[  \rightarrow\left[  0,\mathbb{+\infty}\right[  $
which is continuous, increasing, concave and such that $\omega(0)=0.$

\begin{definition}
Letting $\omega$ be a modulus of continuity, we denote $\sigma\left(
x,\xi\right)  \in\sum_{\omega}$, if $\sigma\left(  x,\xi\right)
:\mathbb{R}^{d}\times\mathbb{R}^{d}\rightarrow\mathbb{C}$ is the continuous
function and for all $\alpha\in\mathbb{N}^{d}$, there exists a constant
$C_{\alpha}$ for which we have
\begin{align}
\left\vert \partial_{\xi}^{\alpha}\sigma\left(  x,\xi\right)  \right\vert  &
\leq A_{\alpha}\left(  1+\left\vert \xi\right\vert \right)  ^{-\left\vert
\alpha\right\vert },\label{eq1.1}\\
\left\vert \partial_{\xi}^{\alpha}\sigma\left(  x+y,\xi\right)  -\partial
_{\xi}^{\alpha}\sigma\left(  x,\xi\right)  \right\vert  &  \leq A_{\alpha
}\omega\left(  \left\vert y\right\vert \right)  \left(  1+\left\vert
\xi\right\vert \right)  ^{-\left\vert \alpha\right\vert }. \label{eq1.2}%
\end{align}

\end{definition}

\begin{definition}
Let $1<p<\infty.$ A measurable function $f$ is said to belong to the weighted
$L^{p},$ $L^{p}\left(  \mathbb{R}^{d},wdx\right)  $, with weight function $w$,
if
\[%
{\displaystyle\int\limits_{\mathbb{R}^{d}}}
\left\vert f(x)\right\vert ^{p}w(x)dx<\infty.
\]
We denote the weighted $L^{p}$ norm by%
\[
\left\Vert f\right\Vert _{L_{w}^{p}}=\left(
{\displaystyle\int\limits_{\mathbb{R}^{d}}}
\left\vert f(x)\right\vert ^{p}w(x)dx\right)  ^{\frac{1}{p}}.
\]

\end{definition}

For $1<p<\infty,$ a positive weight function $w$ is said to be in the class
$A_{p}$ if $w$ is locally integrable and satisfies the condition
\begin{equation}
\underset{Q}{\sup}\left(  \frac{1}{\left\vert Q\right\vert }%
{\displaystyle\int\limits_{Q}}
w(x)dx\right)  \left(  \frac{1}{\left\vert Q\right\vert }%
{\displaystyle\int\limits_{Q}}
w^{-\frac{1}{\left(  p-1\right)  }}(x)dx\right)  ^{p-1}<\infty,\label{eq1.3}%
\end{equation}
where the supremum is taken over all cubes $Q$ in $\mathbb{R}^{d}.$

Our result is stated as follows :

\begin{theorem}
\label{Th1}If the modulus of continuity $\omega$ satisfies the condition
$j^{2}\omega\left(  2^{-j}\right)  <C$, for all $j\in\mathbb{N}$, then any
pseudo-differential operator $\sigma$ with symbol $\sigma\left(  x,\xi\right)
\in\sum_{\omega}$ has a bounded extention to all of $X^{r}=\mathcal{M}\left(
H^{r}\rightarrow L^{2}\right)  $.
\end{theorem}


To prove this theorem, we first introduce some notations. Let $Q$ denote any
cube in $\mathbb{R}^{d}$ and write $\left\vert Q\right\vert $ for the Lebesgue
measure of $Q$. For a locally integrable function $f$, let $f_{Q}$ denote the
mean value of $f$ over $Q$, that is
\[
f_{Q}=\frac{1}{\left\vert Q\right\vert }%
{\displaystyle\int\limits_{Q}}
f(x)dx.
\]

We list the several operators we use later :

\begin{description}
\item[(a)] The Hardy-Littlewood maximal function, $Mf$, for a locally
integrable function $f$ on $\mathbb{R}^{d}$ by
\[
Mf(x)=\underset{Q}{\sup}\frac{1}{\left\vert Q\right\vert }%
{\displaystyle\int\limits_{Q}}
\left\vert f(y)\right\vert dy,
\]
where the supremum ranges over all cubes $Q$ containing $x$.

\item[(b)] Modified maximal function of $f$ :%
\[
M_{\gamma}f(x)=\underset{Q}{\sup}\left(  \frac{1}{\left\vert Q\right\vert }%
{\displaystyle\int\limits_{Q}}
\left\vert f(y)\right\vert ^{\gamma}dy\right)  ^{\frac{1}{\gamma}},
\]
where the supremum is taken over all cubes $Q$ containing $x$.

\item[(c)] Dyadic maximal function of $f$ :%
\[
f^{\ast}(x)=\underset{Q}{\sup}\frac{1}{\left\vert Q\right\vert }%
{\displaystyle\int\limits_{Q}}
\left\vert f(y)\right\vert dy
\]
where the supremum is taken over all dyadic cubes $Q$ , with sides parallel to
the axes containing $x$.

\item[(d)]
\[
f^{\neq}(x)=\underset{Q}{\sup}\frac{1}{\left\vert Q\right\vert }%
{\displaystyle\int\limits_{Q}}
\left\vert f(y)-f_{Q}\right\vert dy,
\]
where the supremum is taken over all cubes $Q$ containing $x$.
\end{description}

\begin{lemma}
[M, Lemma 2.1]\label{Lemma.1}Let $w\in A_{p}$, then $\mathcal{S}\left(
\mathbb{R}^{d}\right)  $ is dense in $L^{p}\left(  \mathbb{R}^{d},wdx\right)
$, $1<p<\infty.$
\end{lemma}

We shall make use of this fact later.

\begin{lemma}
[CM, Theorem 9]Let $\sigma$ be a pseudo-differential operator with symbol
$\sigma\left(  x,\xi\right)  \in\sum_{\omega}$. Then the following two
conditions are equivalent :%
\begin{equation}%
{\displaystyle\sum\limits_{j=0}^{\infty}}
\left[  \omega\left(  2^{-j}\right)  \right]  ^{2}<+\infty, \label{eq1.4}%
\end{equation}%
\begin{equation}
\text{For all }p\text{, }1<p<\infty,\text{ }\sigma\text{ is bounded on }%
L^{p}\left(  \mathbb{R}^{d}\right)  . \label{eq1.5}%
\end{equation}

\end{lemma}

\begin{definition}
[CM, p.41]\label{Def3} We say that $\sigma\left(  x,\xi\right)  \in
\sum_{\omega}$ is a reduced symbol if there exist a constant $C_{1}>0$, a
function $\phi\in C_{0}^{\infty}\left(  \mathbb{R}^{d}\right)  $ and a
sequence $m_{j}$, $j\geq0$ of continuous functions on $\mathbb{R}^{d}$ such
that%
\begin{equation}
\sigma\left(  x,\xi\right)  =%
{\displaystyle\sum\limits_{j=0}^{\infty}}
m_{j}(x)\phi\left(  2^{-j}\xi\right)  , \label{eq1.6}%
\end{equation}
where
\begin{equation}
\left\Vert m_{j}\right\Vert _{L^{\infty}}\leq C_{1} \label{eq1.7}%
\end{equation}%
\begin{equation}
\left\Vert m_{j}(x+y)-m_{j}(x)\right\Vert _{L^{\infty}}\leq C_{1}\omega\left(
\left\vert y\right\vert \right)  \label{eq1.8}%
\end{equation}%
\begin{equation}
\phi\text{ is supported in }\frac{1}{3}\leq\left\vert \xi\right\vert \leq3
\label{eq1.9}%
\end{equation}
and%
\begin{equation}
\left\vert \partial_{\xi}^{\alpha}\phi(\xi)\right\vert \leq C_{1}\text{ for
}\left\vert \alpha\right\vert \leq d. \label{eq1.10}%
\end{equation}

\end{definition}

\begin{lemma}
[CM, Proposition 5, p.46]\label{Lemma.3}For every symbol $\sigma\left(
x,\xi\right)  \in\sum_{\omega}$, we can find a sequence of reduced symbols
$\sigma_{k}\left(  x,\xi\right)  $, $k\in\mathbb{Z}^{d}$, such that%
\[
\sigma\left(  x,\xi\right)  =\tau\left(  x,\xi\right)  +\sum_{k\in
\mathbb{Z}^{d}}\left(  1+\left\vert k\right\vert ^{2}\right)  ^{-d}\sigma
_{k}\left(  x,\xi\right)
\]
and
\begin{align*}
\left\vert \partial_{\xi}^{\alpha}\tau\left(  x,\xi\right)  \right\vert  &
\leq C_{\alpha},\\
\tau\left(  x,\xi\right)   &  =0\text{ \quad if }\left\vert \xi\right\vert
\geq1.
\end{align*}

\end{lemma}

\begin{lemma}
[St1, p.63]\label{Lemma.4}Let $\phi$ be a radial, decreasing, positive and
integrable function. Set $\phi_{t}(x)=t^{-d}\phi\left(  t^{-1}x\right)  $.
Then
\[
\underset{t>0}{\sup}\left\vert \phi_{t}\ast f(x)\right\vert \leq CMf(x)\text{
\quad for \quad}f\in\mathcal{S}\left(  \mathbb{R}^{d}\right)  .
\]

\end{lemma}

\begin{lemma}
[M, Lemma 2.9]\label{Lemma.5}Let $\phi$ be a function in definition
\ref{Def3}. Then for $t\geq0$, there is a constant $C_{t}$ such that
the inequality
\[
\left\vert y\right\vert ^{t}\left\vert
{\displaystyle\int\limits_{\mathbb{R}^{d}}}
\phi\left(  2^{-j}\xi\right)  e^{2\pi ix\xi}d\xi\right\vert \leq
C_{t}2^{j\left(  d-t\right)  }%
\]
holds for all $y\in\mathbb{R}^{d}$ and every integer $j\geq0$.
\end{lemma}

\begin{lemma}
[M, Lemma 2.7]\label{Lemma.6}There is a constant $C>0$ such that
\[
\left\Vert f^{\ast}\right\Vert _{L_{w}^{2}}\leq C\left\Vert f^{\neq
}\right\Vert _{L_{w}^{2}}\text{ for all }f\in L^{2}\left(  \mathbb{R}%
^{d},wdx\right)  \cap L^{1}\left(  \mathbb{R}^{d}\right)  .
\]

\end{lemma}

\section{An $L^{2}\left(  \mathbb{R}^{d},wdx\right)  $ theorem}

After these preliminaries, we state the first main result which constitutes
the main part of the proof of theorem \ref{Th1}.

\begin{theorem}
\label{th2}Suppose $1<\gamma<\infty$ and let $\sigma$ be a
pseudo-differential operator with symbol $\sigma\left(  x,\xi\right)  \in
\sum_{\omega}$. If the modulus of continuity $\omega$ satisfies the condition
\begin{equation}
j^{2}\omega\left(  2^{-j}\right)  <C\text{ \quad for all }j\geq0,
\label{eq1.11}%
\end{equation}
then there is a constant $C>0$ such that the pointwise estimate
\[
\left(  \sigma f\right)  ^{\neq}(x_{0})\leq CM_{\gamma}f(x_{0})
\]
holds for all $x_{0}\in\mathbb{R}^{d}$ and $f\in\mathcal{S}\left(
\mathbb{R}^{d}\right)  .$
\end{theorem}

\begin{proof}
The proof is based on the idea of the proof of theorem 2.8 in [M].

Given $x_{0}\in\mathbb{R}^{d}$, we let $Q$ be a cube containing $x_{0}$, with
center $x^{\prime}$ and diameter $D$. Fix a function $\eta\in C_{0}^{\infty
}\left(  \mathbb{R}^{d}\right)  $ so that $0\leq\eta(x)\leq1$, with
$\eta(x)=1$ for $\left\vert x-x^{\prime}\right\vert \leq2D,$ and $\eta(x)=0$
for $\left\vert x-x^{\prime}\right\vert \geq3D.$ Then, for $f\in
\mathcal{S}\left(  \mathbb{R}^{d}\right)  $
\[
\frac{1}{\left\vert Q\right\vert }%
{\displaystyle\int\limits_{Q}}
\left\vert \sigma f(x)-\left(  \sigma f\right)  _{Q}\right\vert dx\leq\frac
{2}{\left\vert Q\right\vert }%
{\displaystyle\int\limits_{Q}}
\left\vert \sigma\left(  \eta f\right)  (x)\right\vert dx
\]%
\[
+\frac{1}{\left\vert Q\right\vert }%
{\displaystyle\int\limits_{Q}}
\left\vert \sigma\left(  \left(  1-\eta\right)  f\right)  (x)-\left[
\sigma\left(  \left(  1-\eta\right)  f\right)  \right]  _{Q}\right\vert dx.
\]
Letting $Q^{\prime}$ be the cube centered at $x^{\prime}$, with sides parallel
to those of $Q$ and with diameter $4D.$ Since the Hardy-Littlewood maximal
operator is bounded on $L^{\gamma}\left(  \mathbb{R}^{d}\right)  $ for
$1<\gamma<\infty,$ we see that the first term is dominated by
\begin{align*}
\frac{2}{\left\vert Q\right\vert }%
{\displaystyle\int\limits_{Q}}
\left\vert \sigma\left(  \eta f\right)  (x)\right\vert dx  &  \leq2\left(
\frac{1}{\left\vert Q\right\vert }%
{\displaystyle\int\limits_{Q}}
\left\vert \sigma\left(  \eta f\right)  (x)\right\vert ^{\gamma}dx\right)
^{\frac{1}{\gamma}}\\
&  \leq C_{\gamma}\left(  \frac{1}{\left\vert Q\right\vert }%
{\displaystyle\int\limits_{\mathbb{R}^{d}}}
\left\vert \left(  \eta f\right)  (x)\right\vert ^{\gamma}dx\right)
^{\frac{1}{\gamma}}\\
&  \leq C_{\gamma}\left(  \frac{1}{\left\vert Q^{\prime}\right\vert }%
{\displaystyle\int\limits_{Q^{\prime}}}
\left\vert f(x)\right\vert ^{\gamma}dx\right)  ^{\frac{1}{\gamma}}\leq
C_{\gamma}M_{\gamma}f(x_{0}).
\end{align*}

To deal with the second term, we for simplicity write $f$ for $\left(
1-\eta\right)  f$, and we assume that $f$ has the support in the set $\left\{
x:\left\vert x-x^{\prime}\right\vert \geq2D\right\}  .$

We begin by decomposing the symbol $\sigma\left(  x,\xi\right)  $ into the sum
of simpler symbols by making use of Lemma \ref{Lemma.3}. Then we can write%
\begin{align*}
\left(  \sigma f\right)  (x)  &  =%
{\displaystyle\int\limits_{\mathbb{R}^{d}}}
\sigma\left(  x,\xi\right)  \widehat{f}(\xi)e^{2\pi ix\xi}d\xi\\
&  =%
{\displaystyle\int\limits_{\mathbb{R}^{d}}}
\tau\left(  x,\xi\right)  \widehat{f}(\xi)e^{2\pi ix\xi}d\xi\\
&  +%
{\displaystyle\int\limits_{\mathbb{R}^{d}}}
f(y)%
{\displaystyle\int\limits_{\mathbb{R}^{d}}}
\sum_{k\in\mathbb{Z}^{d}}\left(  1+\left\vert k\right\vert ^{2}\right)
^{-d}\sigma_{k}\left(  x,\xi\right)  e^{2\pi i\left(  x-y\right)  \xi}d\xi
dy\\
&  =Tf(x)+\sum_{k\in\mathbb{Z}^{d}}\left(  1+\left\vert k\right\vert
^{2}\right)  ^{-d}S_{k}f(x).
\end{align*}

$T$ is a pseudo-differential operator whose symbol is $\tau\left(
x,\xi\right)  ;$ the $\xi-$support of this symbol is contained in the set
$\left\{  \xi:\left\vert \xi\right\vert \leq1\right\}  ,$ and $\tau\left(
x,\xi\right)  $ has the property that
\[
\left\vert \partial_{\xi}^{\alpha}\tau\left(  x,\xi\right)  \right\vert \leq
C_{\alpha}\text{, \quad}\alpha\in\mathbb{N}^{d}.
\]
Thus, we can write
\[
Tf(x)=%
{\displaystyle\int\limits_{\mathbb{R}^{d}}}
f(y)K\left(  x,x-y\right)  dy,
\]
where we have set
\[
K\left(  x,y\right)  =K\left(  x,x-y\right)  =%
{\displaystyle\int\limits_{\mathbb{R}^{d}}}
\tau\left(  x,\xi\right)  e^{2\pi ix\xi}d\xi.
\]
Note that $K\left(  x,y\right)  $ has the property that
\[
\left\vert K\left(  x,y\right)  \right\vert \leq C_{m}\left(  1+\left\vert
y\right\vert \right)  ^{-m}\text{ \quad for all }x\in\mathbb{R}^{d},
\]
where $m$ is any integer greater than $d$, and $C_{m}$ is a constant
independent of $x$. In fact more generally we have
\begin{align*}
\left\vert y^{\alpha}\partial_{y}^{\beta}K\left(  x,y\right)  \right\vert  &
=A_{\alpha,\beta}\left\vert
{\displaystyle\int\limits_{\mathbb{R}^{d}}}
\tau\left(  x,\xi\right)  \xi^{\beta}\partial_{\xi}^{\alpha}e^{2\pi iy\xi}%
d\xi\right\vert \\
&  \leq A_{\alpha,\beta}%
{\displaystyle\int\limits_{\mathbb{R}^{d}}}
\left\vert \partial_{\xi}^{\alpha}\left[  \tau\left(  x,\xi\right)  \xi
^{\beta}\right]  \right\vert d\xi\\
&  \leq A_{\alpha,\beta},
\end{align*}
with $A_{\alpha,\beta}$ independent of $x$ and $y$. Then by lemma
\ref{Lemma.4}, we have%
\begin{align*}
\left\vert Tf(x)\right\vert  &  \leq%
{\displaystyle\int\limits_{\mathbb{R}^{d}}}
\left\vert f(y)\right\vert \left\vert K\left(  x,x-y\right)  \right\vert dy\\
&  \leq C_{m}%
{\displaystyle\int\limits_{\mathbb{R}^{d}}}
\left\vert f(y)\right\vert \left(  1+\left\vert x-y\right\vert \right)
^{-m}dy\\
&  \leq C_{m}Mf(x).
\end{align*}
Thus, we have
\[
\left(  Tf\right)  ^{\neq}(x_{0})\leq C_{m}M_{\gamma}f(x),
\]
and hence
\[
\left(  \sigma f\right)  ^{\neq}(x_{0})\leq CM_{\gamma}f(x_{0})+\sum
_{k\in\mathbb{Z}^{d}}\left(  1+\left\vert k\right\vert ^{2}\right)
^{-d}\left(  S_{k}f\right)  ^{\neq}(x_{0}).
\]

Therefore, our next task is to examine the operator $S_{k}.$ We note that
$\sigma_{k}\left(  x,\xi\right)  $ satisfies the condition (\ref{eq1.6}) to
(\ref{eq1.10}) in definition \ref{Def3} with $m_{j,k}$, $\phi_{k}$ in
place of $m_{j}$, $\phi$ respectively, where $C^{\prime}$s are independent of
$k$ and also of $j$. Then for every $k$,
\begin{align*}
S_{k}f(x)  &  =%
{\displaystyle\int\limits_{\mathbb{R}^{d}}}
f(y)%
{\displaystyle\int\limits_{\mathbb{R}^{d}}}
\sum_{j=0}^{\infty}m_{j,k}(x)\phi_{k}\left(  2^{-j}\xi\right)  e^{2\pi
i\left(  x-y\right)  \xi}d\xi dy\\
&  =\sum_{j=0}^{\infty}A_{j,k}f(x).
\end{align*}
We now estimate $\left(  S_{k}f\right)  ^{\neq}(x_{0}).$
\begin{equation}
\frac{1}{\left\vert Q\right\vert }%
{\displaystyle\int\limits_{Q}}
\left\vert A_{j,k}f(x)-\left(  A_{j,k}f\right)  _{Q}\right\vert dx=\frac
{1}{\left\vert Q\right\vert }%
{\displaystyle\int\limits_{Q}}
\left\vert \frac{1}{\left\vert Q\right\vert }%
{\displaystyle\int\limits_{Q}}
\left[  A_{j,k}f(x)-\left(  A_{j,k}f\right)  (z)\right]  dz\right\vert
dx\nonumber
\end{equation}%
\begin{equation}
=\frac{1}{\left\vert Q\right\vert }%
{\displaystyle\int\limits_{Q}}
\frac{1}{\left\vert Q\right\vert }\left\vert
{\displaystyle\int\limits_{\mathbb{R}^{d}}}
f(y)%
{\displaystyle\int\limits_{\mathbb{R}^{d}}}
\phi_{k}\left(  2^{-j}\xi\right)  \left[  m_{j,k}(x)e^{2\pi i\left(
x-y\right)  \xi}-m_{j,k}(z)e^{2\pi i\left(  z-y\right)  \xi}\right]  d\xi
dydz\right\vert dx. \label{eq1.12}%
\end{equation}
To estimate this quantity, we consider two cases :

Case 1. $2^{j}D\geq1$. The last quantity is dominated by
\begin{align*}
&  2%
{\displaystyle\sum\limits_{h=1}^{\infty}}
\frac{1}{\left\vert Q\right\vert }%
{\displaystyle\int\limits_{Q}}
{\displaystyle\int\limits_{2^{h}D\leq\left\vert y-x^{\prime}\right\vert
\leq2^{h+1}D}}
\left\vert f(y)\right\vert \left\vert
{\displaystyle\int\limits_{\mathbb{R}^{d}}}
\phi_{k}\left(  2^{-j}\xi\right)  m_{j,k}(x)e^{2\pi i\left(  x-y\right)  \xi
}d\xi\right\vert dydx\\
&  \leq C%
{\displaystyle\sum\limits_{h=1}^{\infty}}
{\displaystyle\int\limits_{Q}}
\frac{2^{hd}}{\left\vert Q_{h}\right\vert }%
{\displaystyle\int\limits_{2^{h}D\leq\left\vert y-x^{\prime}\right\vert
\leq2^{h+1}D}}
\frac{\left\vert f(y)\right\vert }{\left\vert x-y\right\vert ^{d+1}}\left\vert
x-y\right\vert ^{d+1}.\\
&  .\left\vert
{\displaystyle\int\limits_{\mathbb{R}^{d}}}
\phi_{k}\left(  2^{-j}\xi\right)  e^{2\pi i\left(  x-y\right)  \xi}%
d\xi\right\vert \left\vert m_{j,k}(x)\right\vert dydx,
\end{align*}
where $Q_{h}$ is the cube with center $x^{\prime}$ with sides parallel to
those of $Q$ and with diameter $2^{h+2}D.$ The last term is bounded by%
\[
C%
{\displaystyle\sum\limits_{h=1}^{\infty}}
D^{d}2^{hd}\left(  2^{h}D\right)  ^{-d-1}2^{-j}\frac{1}{\left\vert
Q_{h}\right\vert }%
{\displaystyle\int\limits_{Q_{h}}}
\left\vert f(y)\right\vert dy
\]
by lemma \ref{Lemma.5} with $t=d+1,$ and the condition (\ref{eq1.7}) of
$m_{j,k}$%
\[
C%
{\displaystyle\sum\limits_{h=1}^{\infty}}
D^{d}2^{hd}\left(  2^{h}D\right)  ^{-d-1}2^{-j}\frac{1}{\left\vert
Q_{h}\right\vert }%
{\displaystyle\int\limits_{Q_{h}}}
\left\vert f(y)\right\vert dy\leq C\left(  2^{j}D\right)  ^{-1}Mf(x_{0}).
\]

Case 2. $2^{j}D<1.$ In this case, (\ref{eq1.12}) is dominated by%
\begin{align*}
&  \frac{1}{\left\vert Q\right\vert }%
{\displaystyle\int\limits_{Q}}
\frac{1}{\left\vert Q\right\vert }%
{\displaystyle\int\limits_{Q}}
{\displaystyle\sum\limits_{h=1}^{\infty}}
{\displaystyle\int\limits_{2^{h}D\leq\left\vert y-x^{\prime}\right\vert
\leq2^{h+1}D}}
\left\vert f(y)\right\vert .\\
&  .\left\vert
{\displaystyle\int\limits_{\mathbb{R}^{d}}}
\phi_{k}\left(  2^{-j}\xi\right)  \left[  m_{j,k}(x)e^{2\pi i\left(
x-y\right)  \xi}-m_{j,k}(z)e^{2\pi i\left(  z-y\right)  \xi}+m_{j,k}(z)e^{2\pi
i\left(  z-y\right)  \xi}-m_{j,k}(z)e^{2\pi i\left(  z-y\right)  \xi}\right]
d\xi\right\vert dydzdx\\
&  \leq\frac{1}{\left\vert Q\right\vert }%
{\displaystyle\int\limits_{Q}}
\frac{1}{\left\vert Q\right\vert }%
{\displaystyle\int\limits_{Q}}
{\displaystyle\sum\limits_{h=1}^{\infty}}
{\displaystyle\int\limits_{2^{h}D\leq\left\vert y-x^{\prime}\right\vert
\leq2^{h+1}D}}
\left\vert f(y)\right\vert .\\
&  .\left\vert
{\displaystyle\int\limits_{\mathbb{R}^{d}}}
\phi_{k}\left(  2^{-j}\xi\right)  \left[  e^{2\pi i\left(  x-y\right)  \xi
}-e^{2\pi i\left(  z-y\right)  \xi}\right]  m_{j,k}(x)d\xi\right\vert
dydzdx+\frac{1}{\left\vert Q\right\vert }%
{\displaystyle\int\limits_{Q}}
\frac{1}{\left\vert Q\right\vert }%
{\displaystyle\int\limits_{Q}}
{\displaystyle\sum\limits_{h=1}^{\infty}}
{\displaystyle\int\limits_{2^{h}D\leq\left\vert y-x^{\prime}\right\vert
\leq2^{h+1}D}}
\left\vert f(y)\right\vert .\\
&  .\left\vert
{\displaystyle\int\limits_{\mathbb{R}^{d}}}
\phi_{k}\left(  2^{-j}\xi\right)  e^{2\pi i\left(  z-y\right)  \xi}\left[
m_{j,k}(x)-m_{j,k}(z)\right]  d\xi\right\vert dydzdx\\
&  =A+B.
\end{align*}
We first estimate $A$.%
\begin{align*}
A  &  \leq\frac{1}{\left\vert Q\right\vert }%
{\displaystyle\int\limits_{Q}}
\frac{1}{\left\vert Q\right\vert }%
{\displaystyle\int\limits_{Q}}
{\displaystyle\sum\limits_{h=1}^{\infty}}
{\displaystyle\int\limits_{2^{h}D\leq\left\vert y-x^{\prime}\right\vert
\leq2^{h+1}D}}
\left\vert f(y)\right\vert .\\
&  .\left\vert
{\displaystyle\int\limits_{\mathbb{R}^{d}}}
\phi_{k}\left(  2^{-j}\xi\right)
{\displaystyle\sum\limits_{p=1}^{d}}
\left(  x_{p}-z_{p}\right)
{\displaystyle\int\limits_{0}^{1}}
2\pi i\xi_{p}e^{2\pi i\left(  x(t)-y\right)  \xi}d\xi\right\vert
dydz\left\vert m_{j,k}(x)\right\vert dx
\end{align*}
where $x(t)=z+t(x-z)$%
\begin{align*}
&  \leq\frac{1}{\left\vert Q\right\vert }%
{\displaystyle\int\limits_{Q}}
\frac{1}{\left\vert Q\right\vert }%
{\displaystyle\int\limits_{Q}}
{\displaystyle\sum\limits_{h=1}^{\infty}}
{\displaystyle\int\limits_{2^{h}D\leq\left\vert y-x^{\prime}\right\vert
\leq2^{h+1}D}}
\frac{\left\vert f(y)\right\vert }{\left\vert y-x^{\prime}\right\vert
^{d+1/2}}%
{\displaystyle\sum\limits_{p=1}^{d}}
\left\vert x_{p}-z_{p}\right\vert
{\displaystyle\int\limits_{0}^{1}}
\left\vert x(t)-y\right\vert ^{d+1/2}.\\
&  .\left\vert
{\displaystyle\int\limits_{\mathbb{R}^{d}}}
\phi_{k}\left(  2^{-j}\xi\right)  2\pi i\xi_{p}e^{2\pi i\left(  x(t)-y\right)
\xi}d\xi\right\vert dtdydz\left\vert m_{j,k}(x)\right\vert dx.
\end{align*}
The integral with respect to $\xi$ is handled just as in the proof of Lemma
\ref{Lemma.5} with $t=d+\frac{1}{2}$ and we see that the last member is not
greater than
\begin{align*}
&  C%
{\displaystyle\sum\limits_{h=1}^{\infty}}
\left(  2^{h}D\right)  ^{d}\left(  2^{h}D\right)  ^{-d-1/2}D2^{\frac{j}{2}%
}\frac{1}{\left\vert Q_{h}\right\vert }%
{\displaystyle\int\limits_{Q_{h}}}
\left\vert f(y)\right\vert dy\\
&  \leq C\left(  2^{j}D\right)  ^{\frac{1}{2}}Mf(x_{0}).
\end{align*}

Next, we estimate $B$.%
\begin{align*}
B  &  \leq\frac{1}{\left\vert Q\right\vert }%
{\displaystyle\int\limits_{Q}}
\frac{1}{\left\vert Q\right\vert }%
{\displaystyle\int\limits_{Q}}
{\displaystyle\sum\limits_{h=1}^{\infty}}
{\displaystyle\int\limits_{2^{h}D\leq\left\vert y-x^{\prime}\right\vert
\leq2^{h+1}D}}
\left\vert f(y)\right\vert \left\vert
{\displaystyle\int\limits_{\mathbb{R}^{d}}}
\phi_{k}\left(  2^{-j}\xi\right)  e^{2\pi i\left(  z-y\right)  \xi}%
d\xi\right\vert dy\left[  m_{j,k}(x)-m_{j,k}(z)\right]  dzdx\\
&  \leq\frac{1}{\left\vert Q\right\vert }%
{\displaystyle\int\limits_{Q}}
\frac{1}{\left\vert Q\right\vert }%
{\displaystyle\int\limits_{Q}}
{\displaystyle\sum\limits_{h=1}^{N}}
{\displaystyle\int\limits_{2^{h}D\leq\left\vert y-x^{\prime}\right\vert
\leq2^{h+1}D}}
\left\vert f(y)\right\vert \left\vert
{\displaystyle\int\limits_{\mathbb{R}^{d}}}
\phi_{k}\left(  2^{-j}\xi\right)  e^{2\pi i\left(  z-y\right)  \xi}%
d\xi\right\vert dy\left[  m_{j,k}(x)-m_{j,k}(z)\right]  dzdx\\
&  +\frac{1}{\left\vert Q\right\vert }%
{\displaystyle\int\limits_{Q}}
\frac{1}{\left\vert Q\right\vert }%
{\displaystyle\int\limits_{Q}}
{\displaystyle\sum\limits_{h=N+1}^{\infty}}
{\displaystyle\int\limits_{2^{h}D\leq\left\vert y-x^{\prime}\right\vert
\leq2^{h+1}D}}
\left\vert f(y)\right\vert \left\vert
{\displaystyle\int\limits_{\mathbb{R}^{d}}}
\phi_{k}\left(  2^{-j}\xi\right)  e^{2\pi i\left(  z-y\right)  \xi}%
d\xi\right\vert dy\left[  m_{j,k}(x)-m_{j,k}(z)\right]  dzdx\\
&  =B_{1}+B_{2},
\end{align*}
where $N$ is the integer which satisfies $2^{N}D<1\leq2^{N+1}D.$
\begin{align*}
B_{1}  &  =\frac{1}{\left\vert Q\right\vert }%
{\displaystyle\int\limits_{Q}}
\frac{1}{\left\vert Q\right\vert }%
{\displaystyle\int\limits_{Q}}
{\displaystyle\sum\limits_{h=1}^{N}}
{\displaystyle\int\limits_{2^{h}D\leq\left\vert y-x^{\prime}\right\vert
\leq2^{h+1}D}}
\frac{\left\vert f(y)\right\vert }{\left\vert z-y\right\vert ^{d}}\left\vert
z-y\right\vert ^{d}.\\
&  .\left\vert
{\displaystyle\int\limits_{\mathbb{R}^{d}}}
\phi_{k}\left(  2^{-j}\xi\right)  e^{2\pi i\left(  z-y\right)  \xi}%
d\xi\right\vert dy\left[  m_{j,k}(x)-m_{j,k}(z)\right]  dzdx\\
&  \leq C^{\prime}%
{\displaystyle\sum\limits_{h=1}^{N}}
\left(  2^{h}D\right)  ^{d}\left(  2^{h}D\right)  ^{-d}\omega\left(  D\right)
\frac{1}{\left\vert Q_{h}\right\vert }%
{\displaystyle\int\limits_{Q_{h}}}
\left\vert f(y)\right\vert dy
\end{align*}
by lemma \ref{Lemma.5} with $t=d$ and the condition (\ref{eq1.8}) of $m_{j,k}$%
\[
B_{1}\leq CN\omega\left(  2^{-N}\right)  Mf(x_{0}).
\]

\begin{align*}
B_{2}  &  =\frac{1}{\left\vert Q\right\vert }%
{\displaystyle\int\limits_{Q}}
\frac{1}{\left\vert Q\right\vert }%
{\displaystyle\int\limits_{Q}}
{\displaystyle\sum\limits_{h=N+1}^{\infty}}
{\displaystyle\int\limits_{2^{h}D\leq\left\vert y-x^{\prime}\right\vert
\leq2^{h+1}D}}
\frac{\left\vert f(y)\right\vert }{\left\vert z-y\right\vert ^{d}}\left\vert
z-y\right\vert ^{d}.\\
&  .\left\vert
{\displaystyle\int\limits_{\mathbb{R}^{d}}}
\phi_{k}\left(  2^{-j}\xi\right)  e^{2\pi i\left(  z-y\right)  \xi}%
d\xi\right\vert dy\left[  m_{j,k}(x)-m_{j,k}(z)\right]  dzdx\\
&  \leq C^{"}%
{\displaystyle\sum\limits_{h=N+1}^{\infty}}
\left(  2^{h}D\right)  ^{d}\left(  2^{h}D\right)  ^{-d-1}2^{-j}\omega\left(
D\right)  \frac{1}{\left\vert Q_{h}\right\vert }%
{\displaystyle\int\limits_{Q_{h}}}
\left\vert f(y)\right\vert dy
\end{align*}
by lemma \ref{Lemma.5} with $t=d+1$ and the condition (\ref{eq1.8}) of
$m_{j,k}$%
\begin{align*}
B_{2}  &  \leq C^{"}%
{\displaystyle\sum\limits_{h=N+1}^{\infty}}
\left(  2^{h}D\right)  ^{-1}2^{-j}\omega\left(  1\right)  Mf(x_{0})\\
&  \leq C^{"}2^{-j}Mf(x_{0}).
\end{align*}
Thus we have
\begin{align*}
B  &  \leq B_{1}+B_{2}\leq C^{\prime}N\omega\left(  2^{-N}\right)
Mf(x_{0})+C^{"}2^{-j}Mf(x_{0})\\
&  \leq\left\{  C^{\prime}N\omega\left(  2^{-N}\right)  +C^{"}2^{-j}\right\}
Mf(x_{0}).
\end{align*}

Putting two cases together, we have shown that if $Q$ is any cube containing
$x_{0}$, then
\begin{align*}
\left(  S_{k}f\right)  ^{\neq}(x_{0})  &  \leq%
{\displaystyle\sum\limits_{j=0}^{\infty}}
\left(  A_{j,k}f\right)  ^{\neq}(x)\\
&  \leq\left\{  C%
{\displaystyle\sum\limits_{2^{j}D\geq1}}
\left(  2^{j}D\right)  ^{-1}+C^{\prime}%
{\displaystyle\sum\limits_{2^{j}D<1}}
\left(  2^{j}D\right)  ^{\frac{1}{2}}+%
{\displaystyle\sum\limits_{2^{j}D<1}}
\left(  C^{\prime}N\omega\left(  2^{-N}\right)  +C^{"}2^{-j}\right)  \right\}
Mf(x_{0})\\
&  \leq\left\{  C+N^{2}\omega\left(  2^{-N}\right)  \right\}  Mf(x_{0})\leq
CMf(x_{0})\leq CM_{\gamma}f(x_{0}).
\end{align*}
We thus find that
\begin{align*}
\left(  \sigma f\right)  ^{\neq}(x_{0})  &  \leq CM_{\gamma}f(x_{0}%
)+\sum_{k\in\mathbb{Z}^{d}}\left(  1+\left\vert k\right\vert ^{2}\right)
^{-d}CM_{\gamma}f(x_{0})\\
&  \leq CM_{\gamma}f(x_{0}).
\end{align*}
Summarizing, we have shown that if $Q$ is any cube containing $x_{0}$, then%
\begin{align*}
&  \frac{1}{\left\vert Q\right\vert }%
{\displaystyle\int\limits_{Q}}
\left\vert \sigma f(x)-\left(  \sigma f\right)  _{Q}\right\vert dx\\
&  \leq\left(  \sigma f\right)  ^{\neq}(x_{0})+\left[  T\left(  \left(
1-\eta\right)  f\right)  \right]  ^{\neq}(x_{0})+\sum_{k\in\mathbb{Z}^{d}%
}\left(  1+\left\vert k\right\vert ^{2}\right)  ^{-d}\left[  S_{k}\left(
\left(  1-\eta\right)  f\right)  \right]  ^{\neq}(x_{0})\\
&  \leq CM_{\gamma}f(x_{0})+CM_{\gamma}\left(  \left(  1-\eta\right)
f\right)  (x_{0})\\
&  \leq CM_{\gamma}f(x_{0}),
\end{align*}
where the constant $C$ is independent of $Q$, $f$ and $x_{0}.$ When we take
the supremum of the left side over all cubes $Q$ containing $x_{0}$, we
finally obtain the desired inequality :%
\[
\left(  \sigma f\right)  ^{\neq}(x_{0})\leq CM_{\gamma}f(x_{0})
\]
for all $x_{0}\in\mathbb{R}^{d}$, $f\in\mathcal{S}\left(  \mathbb{R}%
^{d}\right)  .$
\end{proof}

We are ready to prove a basic result about pseudo-differential operators.

\begin{theorem}
If $w\in A_{2}\left(  \mathbb{R}^{d}\right)  $ and if the modulus of
continuity $\omega$ satisfies the condition $j^{2}\omega\left(  2^{-j}\right)
<C$, for all $j\in\mathbb{N}$, then any pseudo-differential operator $\sigma$
with symbol $\sigma\left(  x,\xi\right)  \in\sum_{\omega}$, initially defined
on $\mathcal{S}$, extends to a bounded operator from $L^{2}\left(
\mathbb{R}^{d},wdx\right)  $ to itself.
\end{theorem}

To prove the theorem, it suffices to show that
\[
\left\Vert \sigma f\right\Vert _{L_{w}^{2}}\leq C\left\Vert f\right\Vert
_{L_{w}^{2}}\text{, \quad whenever }f\in\mathcal{S}\text{,}%
\]
with $C$ independent of $f$.

\begin{proof}
We prove this in the same way as was used by [M, theorem 2.12]. If
$f\in\mathcal{S}\left(  \mathbb{R}^{d}\right)  ,$ then since $\sigma f\in
L^{2}\left(  \mathbb{R}^{d},wdx\right)  \cap L^{1}\left(  \mathbb{R}%
^{d}\right)  $
\begin{align*}
\left\Vert \sigma f\right\Vert _{L_{w}^{2}}  &  \leq\left\Vert \left(  \sigma
f\right)  ^{\ast}\right\Vert _{L_{w}^{2}}\leq C\left\Vert \left(  \sigma
f\right)  ^{\neq}\right\Vert _{L_{w}^{2}}\\
&  \leq C\left\Vert M_{\gamma}f\right\Vert _{L_{w}^{2}}\text{ , if }%
1<\gamma<\infty\\
&  \leq C\left\Vert f\right\Vert _{L_{w}^{2}}\text{ , if }1<\gamma<2.
\end{align*}
The first inequality is easy, since
\[
\left\vert \sigma f(x)\right\vert \leq\left(  \sigma f\right)  ^{\ast
}(x)\text{ \quad for every }x.
\]
Since $\sigma f\in\mathcal{S}$, $\sigma f\in L^{2}\left(  \mathbb{R}%
^{d},wdx\right)  \cap L^{1}\left(  \mathbb{R}^{d}\right)  ;$ so we can apply
lemma \ref{Lemma.6} to prove the second inequality. The third inequality is
theorem \ref{th2}, while the last inequality is proved like this :%
\begin{align*}
\left\Vert M_{\gamma}f\right\Vert _{L_{w}^{2}}  &  =\left\Vert \left[
M\left\vert f\right\vert ^{\gamma}\right]  ^{\frac{1}{\gamma}}\right\Vert
_{L_{w}^{2}}=\left(
{\displaystyle\int}
\left[  M\left(  \left\vert f\right\vert ^{\gamma}\right)  \right]  ^{\frac
{2}{\gamma}}w(x)dx\right)  ^{\frac{1}{2}}\\
&  \leq C\left(
{\displaystyle\int}
\left\vert f\right\vert ^{2}w(x)dx\right)  ^{\frac{1}{2}}\text{ since }%
\gamma<2\\
&  =C\left\Vert f\right\Vert _{L_{w}^{2}}.
\end{align*}
Because of lemma \ref{Lemma.1}, we can extend $\sigma$ to a bounded operator
on $L^{2}\left(  \mathbb{R}^{d},wdx\right)  .$
\end{proof}

\section{Pointwise multipliers $X^{r}$}

In this section, we give a description of the multiplier space\ $X^{r}$
introduced recently by P.G.\ Lemari\'{e}-Rieusset in his work \cite{Lem}. The
space $X^{r}$of pointwise multipliers which map $L^{2}$ into $H^{-r}$ is
defined in the following way

\begin{definition}
For $0\leq r<\frac{d}{2}$, we define the space $X^{r}\left(  \mathbb{R}%
^{d}\right)  $ as the space of functions, which are locally square integrable
on $\mathbb{R}^{d}$ and such that pointwise multiplication with these
functions maps boundedly $H^{r}\left(  \mathbb{R}^{d}\right)  $ to
$L^{2}\left(  \mathbb{R}^{d}\right)  $, i.e.,
\[
X^{r}=\left\{  f\in L_{loc}^{2}:\text{ }\forall g\in H^{r}\text{ \ }fg\in
L^{2}\right\}
\]
where we denote by $H^{r}\left(  \mathbb{R}^{d}\right)  $ the completion of
the space $\mathcal{D}\left(  \mathbb{R}^{d}\right)  $ with respect to the
norm $\left\Vert u\right\Vert _{H^{r}}=\left\Vert \left(  1-\Delta\right)
^{\frac{r}{2}}u\right\Vert _{L^{2}}.$
\end{definition}

The norm of $X^{r}$ is given by the operator norm of pointwise multiplication
:
\[
\left\Vert f\right\Vert _{X^{r}}=\underset{\left\Vert g\right\Vert _{H^{r}%
}\leq1}{\sup}\left\Vert fg\right\Vert _{L^{2}}%
\]
We now turn to another way of introducing capacity.

\begin{definition}
[Capacitary measures and capacitary potentials]The Bessel capacity $cap\left(
e;H^{r}\right)  $ of a compact set $e\subset$ $\mathbb{R}^{d}$ is defined by
[AH]
\[
cap\left(  e;H^{r}\right)  =\inf\left\{  \left\Vert u\right\Vert _{H^{r}}%
^{2}\text{ : }u\in C_{0}^{\infty}\left(  \mathbb{R}^{d}\right)  \text{, }%
u\geq1\text{ sur }e\right\}
\]

\end{definition}


We shall show the following theorem.

\begin{theorem}
\label{th1.1}There is a positive constant $C$ depending only on $d$ such that
\[%
{\displaystyle\int\limits_{0}^{\infty}}
cap\left(  A_{t},H^{r}\right)  d\left(  t^{2}\right)  \leq C\left\Vert
u\right\Vert _{H^{r}}^{2}%
\]
where
\[
A_{t}=\left\{  x\in\mathbb{R}^{d}:u(x)\geq t\right\}  \text{ et }u=G_{r}\ast
f\text{ / }f\in C_{0}^{\infty}\left(  \mathbb{R}^{d}\right)
\]
for any nonnegative measurable function $f$.
\end{theorem}

This theorem was established first by Hansson \cite{Han}. Later Maz'ya
(\cite{Maz}, th.8.2.3) and Adams (\cite{Ad}, th 1.6). Maz'ya and Adams used
the joint measurability of $G_{r}\ast\mu_{t}$ on $\mathbb{R}^{d}$ where
$\mu_{t}$ is the capacitary measure for the set $\left\{  x\in\mathbb{R}%
^{d}:u(x)\geq t\right\}  .$ However, the measurability does not seem to be
obvious. We shall give an elementary proof which gets around this difficulty.

An easy corollary to theorem \ref{th1.1}, we obtain the following
characterization of Carleson types measures.

\begin{corollary}
\label{corol1.2.20}For a nonnegative measure $\mu$, lthe following assertions
are equivalent :

\begin{enumerate}
\item For any $f\in$ $L_{+}^{2}\left(  \mathbb{R}^{d}\right)  ,$ we have
\[
\int\left(  G_{r}\ast f\right)  ^{2}d\mu\leq C_{1}\int f^{2}dx
\]

\item For any compact set $e\subset\mathbb{R}^{d}$, we have
\[
\mu\left(  e\right)  \leq C_{2}cap\left(  e\right)
\]

\end{enumerate}
\end{corollary}

\bigskip Moreover, we have the following characterization :
\[
\int u^{2}d\mu\leq C\left\Vert u\right\Vert _{H^{r}}^{2}\text{ }%
\Leftrightarrow\mu\left(  e\right)  \leq C_{2}cap\left(  e\right)
\]

\begin{remark}
Let $\mu$ nonnegative measure. The inequality%
\[
\int_{\mathbb{R}^{d}}u^{2}d\mu\leq C\left\Vert u\right\Vert _{H^{r}}^{2}\text{
}%
\]
for $u\in C_{0}^{\infty}\left(  \mathbb{R}^{d}\right)  $ is called the trace inequality.
\end{remark}

Before to prove this theorem \ref{th1.1}, we prove the corollary
\ref{corol1.2.20} (and thus the theorem \ref{th1.1}) gives the
characterization of the multipliers spaces.

\begin{proof}
(1)$\Rightarrow$(2). This part can be proved without the capacity inequality.
Let $e$ be a compact set. Take $f\in$ $L_{+}^{2}\left(  \mathbb{R}^{d}\right)
$ such that $G_{r}\ast f$ $\geq1$ on $e$. Then
\[
\mu\left(  e\right)  \leq\int\left(  G_{r}\ast f\right)  ^{2}d\mu\leq
C_{1}\int f^{2}dx
\]
Taking the infinimum with respect to $f$, we obtain
\[
\mu\left(  e\right)  \leq C_{1}cap\left(  e\right)
\]

(2)$\Rightarrow$(1). By the capacitability, we have
\[
\mu\left(  e\right)  \leq C_{2}cap\left(  e\right)
\]
for every Borel set $e$. Let $f\in$ $L_{+}^{2}\left(  \mathbb{R}^{d}\right)  $
and apply the above inequality to $A_{t}:$%
\[
A_{t}=\left\{  x:\text{ }G_{r}\ast f\text{ }(x)\geq t\right\}  .
\]
By theorem \ref{th1.1}, we have%
\begin{align*}
\int f^{2}dx  &  \geq\frac{1}{c}%
{\displaystyle\int\limits_{0}^{\infty}}
cap\left(  A_{t},H^{r}\right)  d\left(  t^{2}\right)  \geq\frac{1}{c}%
{\displaystyle\int\limits_{0}^{\infty}}
\mu\left(  A_{t}\right)  d\left(  t^{2}\right) \\
&  =\frac{1}{c}\int\left(  G_{r}\ast f\right)  ^{2}d\mu
\end{align*}

\end{proof}

To proof theorem \ref{th1.1}, we will need several lemmas.

\begin{lemma}
\label{lem1.1}If $a_{j}\geq0$ for $j\in\mathbb{Z}$, then
\[
\left(  \sum_{\text{ }j\in\mathbb{Z}}a_{j}\right)  ^{2}\leq2\sum_{\text{ }%
i\in\mathbb{Z}}a_{i}\left(  \sum_{\text{ }j=i}^{\infty}a_{j}\right)
\]

\end{lemma}

The proof is immediat.

\begin{lemma}
\label{lem1.2}Suppose $\mu_{j}$ are measures function such that
\[
G_{r}\ast\left(  G_{r}\ast\mu_{j}\right)  \leq1\text{ on supp}\left(  \mu
_{j}\right)  \text{ for }j\in\mathbb{Z}%
\]
Then
\[%
{\displaystyle\int\limits_{\mathbb{R}^{d}}}
\sum_{\text{ }j\in\mathbb{Z}}\left(  2^{j}G_{r}\ast\mu_{j}\right)  ^{2}dx\leq
c\sum_{\text{ }j\in\mathbb{Z}}2^{2j}\left\Vert \mu_{j}\right\Vert
\]

\end{lemma}

\begin{proof}
Apply the equilibrium potential of $\mu_{j}$
\[
G_{r}\ast\left(  G_{r}\ast\mu_{j}\right)  \leq K\text{ on }\mathbb{R}^{d}%
\]
to obtain%
\begin{align*}%
{\displaystyle\int\limits_{\mathbb{R}^{d}}}
\sum_{\text{ }j\in\mathbb{Z}}\left(  2^{j}G_{r}\ast\mu_{j}\right)  ^{2}dx  &
\leq2%
{\displaystyle\int\limits_{\mathbb{R}^{d}}}
\sum_{i\in\mathbb{Z}}\left(  2^{i}G_{r}\ast\mu_{i}\right)  \sum_{j=i}%
^{+\infty}\left(  2^{j}G_{r}\ast\mu_{j}\right)  dx\\
&  =2\sum_{\text{ }j\in\mathbb{Z}}2^{j}\sum_{i=-\infty}^{j}2^{i}%
{\displaystyle\int\limits_{\mathbb{R}^{d}}}
\left(  G_{r}\ast G_{r}\ast\mu_{i}\right)  d\mu_{j}\\
&  \leq2\sum_{\text{ }j\in\mathbb{Z}}2^{j}\sum_{i=-\infty}^{j}2^{i}M%
{\displaystyle\int\limits_{\mathbb{R}\ ^{d}}}
d\mu_{j}\\
&  =c\sum_{\text{ }j\in\mathbb{Z}}2^{2j}\left\Vert \mu_{j}\right\Vert .
\end{align*}
Hence,%
\[%
{\displaystyle\int\limits_{\mathbb{R}^{d}}}
\sum_{\text{ }j\in\mathbb{Z}}\left(  2^{j}G_{r}\ast\mu_{j}\right)  ^{2}dx\leq
c\sum_{\text{ }j\in\mathbb{Z}}2^{2j}\left\Vert \mu_{j}\right\Vert .
\]
The lemma follows.
\end{proof}

\begin{lemma}
\bigskip Let $f$ be a nonnegative continuous function of compact support and
$e$ a Borel set. \ Let
\[
A_{j}=\left\{  x\in e:u(x)\geq2^{j}\right\}
\]
and let $\mu_{j}$ be the capacitary measure for $A_{j}$, i.e.,

\begin{enumerate}
\item supp$\left(  \mu_{j}\right)  \subset\overline{A_{j}}$,

\item $G_{r}\ast\left(  G_{r}\ast\mu_{i}\right)  \geq1$, p.p. sur $A_{j}$

\item $G_{r}\ast\left(  G_{r}\ast\mu_{i}\right)  \leq1,$ sur supp$\left(
\mu_{j}\right)  $

\item $\left\Vert \mu_{j}\right\Vert =cap\left(  A_{j}\right)  $.
\end{enumerate}
\end{lemma}

Then%
\[
\frac{3}{4}\sum_{\text{ }j\in\mathbb{Z}}2^{2j}\left\Vert \mu_{j}\right\Vert
\leq%
{\displaystyle\int\limits_{0}^{\infty}}
cap\left(  \left\{  x:u(x)\geq t\right\}  ;H^{r}\right)  d\left(
t^{2}\right)  \leq3\left\Vert f\right\Vert _{L^{2}}\left\Vert \sum_{\text{
}j\in\mathbb{Z}}2^{j}\left(  G_{r}\ast\mu_{j}\right)  \right\Vert _{L^{2}}%
\]

\begin{proof}
By definition%
\[%
{\displaystyle\int\limits_{0}^{\infty}}
cap\left(  \left\{  x:u(x)\geq t\right\}  ;H^{r}\right)  d\left(
t^{2}\right)  =\sum_{\text{ }j\in\mathbb{Z}}%
{\displaystyle\int\limits_{2^{j}}^{2^{j+1}}}
cap\left(  \left\{  x:u(x)\geq t\right\}  ;H^{r}\right)  d(t^{2})
\]
Or%
\[
\left(  2^{2\left(  j+1\right)  }-2^{2j}\right)  cap\left(  A_{j+1}\right)
\leq%
{\displaystyle\int\limits_{2^{j}}^{2^{j+1}}}
cap\left(  \left\{  x:u(x)\geq t\right\}  ;H^{r}\right)  d(t^{2})\leq\left(
2^{2\left(  j+1\right)  }-2^{2j}\right)  cap\left(  A_{j}\right)
\]
The left hand side is equal to%
\[
2^{2\left(  j+1\right)  }\left(  1-\frac{1}{4}\right)  cap\left(
A_{j+1}\right)  =\frac{3}{4}2^{2\left(  j+1\right)  }cap\left(  A_{j+1}%
\right)  =\frac{3}{4}2^{2\left(  j+1\right)  }\left\Vert \mu_{j+1}\right\Vert
\]
which implies the first required inequality, i.e.,%
\[
\frac{3}{4}\sum_{\text{ }j\in\mathbb{Z}}2^{2j+2}\left\Vert \mu_{j+1}%
\right\Vert \leq\sum_{\text{ }j\in\mathbb{Z}}%
{\displaystyle\int\limits_{2^{j}}^{2^{j+1}}}
cap\left(  \left\{  x:u(x)\geq t\right\}  ;H^{r}\right)  d(t^{2})
\]
The right hand side is equal to%
\begin{align*}
\left(  2^{2\left(  j+1\right)  }-2^{2j}\right)  cap\left(  A_{j}\right)   &
=3\times2^{2j}\left\Vert \mu_{j}\right\Vert \leq3\times2^{2j}\int
_{\mathbb{R}^{d}}\frac{G_{r}\ast f}{2^{j}}d\mu_{j}\\
&  =3\times2^{j}\int_{\mathbb{R}^{d}}f\left(  G_{r}\ast\mu_{j}\right)  dx
\end{align*}
Adding this, we obtain%
\begin{align*}%
{\displaystyle\int\limits_{0}^{\infty}}
cap\left(  \left\{  x:u(x)\geq t\right\}  ;H^{r}\right)  d\left(
t^{2}\right)   &  \leq3\sum_{\text{ }j\in\mathbb{Z}}2^{j}\int_{\mathbb{R}^{d}%
}f\left(  G_{r}\ast\mu_{j}\right)  dx\\
&  =3\int_{\mathbb{R}^{d}}f\left(  \sum_{\text{ }j\in\mathbb{Z}}2^{j}\left(
G_{r}\ast\mu_{j}\right)  \right)  dx
\end{align*}
Hence H\"{o}lder's inequality yields%
\[%
{\displaystyle\int\limits_{0}^{\infty}}
cap\left(  \left\{  x:u(x)\geq t\right\}  ;H^{r}\right)  d\left(
t^{2}\right)  \leq3\left\Vert f\right\Vert _{L^{2}}\left\Vert \sum_{\text{
}j\in\mathbb{Z}}2^{j}\left(  G_{r}\ast\mu_{j}\right)  \right\Vert _{L^{2}}%
\]

\end{proof}

\bigskip We now are in position to prove the theorem \ref{th1.1}.
For the preceding corollary, we are immediately deduce the main result.

\begin{proof}
In view of the monotone convergence theorem, it is sufficient to show that%
\[
I=%
{\displaystyle\int\limits_{0}^{\infty}}
cap\left(  \left\{  x\in B\left(  0,R\right)  :u(x)\geq t\right\}  \right)
d\left(  t^{2}\right)  \leq C\left\Vert f\right\Vert _{L^{2}}^{2}%
\]
for $R>0$, $f\in C_{0}^{\infty}\left(  \mathbb{R}^{d}\right)  $ and $C$ is
independent of $R$ and $f$. Since $G_{r}\ast f=u$ is bounded, it follows that
$\left\{  x:u(x)\geq t\right\}  =\emptyset$ for large $t$, say $t>T,$ so that
\[
I\leq cap\left(  B\left(  0,R\right)  \right)  T^{2}<\infty.
\]
Lett $\mu_{j}$ be the capacitary measure for $A_{j}$%
\[
A_{j}=\left\{  x\in B\left(  0,R\right)  :u(x)\geq2^{j}\right\}
\]
By lemma \ref{lem1.2} with $A=B\left(  0,R\right)  $, we have%
\begin{equation}
I\leq3\left\Vert f\right\Vert _{L^{2}}\left\Vert \sum_{\text{ }j\in\mathbb{Z}%
}2^{j}\left(  G_{r}\ast\mu_{j}\right)  \right\Vert _{L^{2}} \label{eq1.8}%
\end{equation}
Then lemma \ref{lem1.2} yield
\[
\left\Vert \sum_{\text{ }j\in\mathbb{Z}}2^{j}\left(  G_{r}\ast\mu_{j}\right)
\right\Vert _{L^{2}}\leq c\left(  \sum_{\text{ }j\in\mathbb{Z}}2^{2j}%
\left\Vert \mu_{j}\right\Vert \right)  ^{\frac{1}{2}}\text{ }\left(
\text{lemma \ref{lem1.2}}\right)
\]
\ and $\ $%
\[
\frac{3}{4}\sum_{\text{ }j\in\mathbb{Z}}2^{2j}\left\Vert \mu_{j}\right\Vert
_{L^{1}}\leq I.
\]
Then%
\[
\left\Vert \sum_{\text{ }j\in\mathbb{Z}}2^{j}\left(  G_{r}\ast\mu_{j}\right)
\right\Vert _{L^{2}}\leq cI^{\frac{1}{2}}%
\]
Since $I<\infty,$ it follows from (\ref{eq1.8}) that
\[
I\leq c\left\Vert f\right\Vert _{L^{2}}I^{\frac{1}{2}}%
\]
i.e
\[
I\leq c\left\Vert f\right\Vert _{L^{2}}^{2}%
\]
Finally,%
\[%
{\displaystyle\int\limits_{0}^{\infty}}
cap\left(  \left\{  x\in B\left(  0,R\right)  :u(x)\geq t\right\}  \right)
d\left(  t^{2}\right)  \leq C\left\Vert f\right\Vert _{L^{2}}^{2}%
\]
The theorem is proved.
\end{proof}

\begin{remark}
Under this conditions and the preceding results, let $\mu
=f^{2}$, it follows that
\[
\int_{\mathbb{R}^{d}}u^{2}f^{2}dx\leq C\left\Vert u\right\Vert _{H^{r}}%
^{2}\text{ }%
\]
and hence
\[%
{\displaystyle\int\limits_{e}}
f^{2}dx\leq C_{2}cap\left(  e,H^{r}\right)  .
\]
\bigskip
\end{remark}

Then, we define the norme $\left\Vert f\right\Vert _{\mathcal{M}\left(
H^{r}\rightarrow L^{2}\right)  }$ by \cite{MS}%
\begin{equation}
\left\Vert f\right\Vert _{\mathcal{M}\left(  H^{r}\rightarrow L^{2}\right)
}\sim\underset{e}{\sup}\frac{\left(
{\displaystyle\int\limits_{e}}
f^{2}dx\right)  ^{\frac{1}{2}}}{\left[  cap\left(  e\text{,}H^{r}\right)
\right]  ^{\frac{1}{2}}}\text{ } \label{eq1.11}%
\end{equation}

We will need the following theorem, which shows that many operators of
classical analysis are bounded in the space of multipliers.

\begin{theorem}
\label{thMV}Let $0\leq r<\frac{d}{2}.$ Suppose that a function\ $h\in
L_{loc\text{ }}^{2}$satisfies%
\begin{equation}
\int\limits_{e}\left\vert h(x)\right\vert ^{2}dx\leq Ccap(e) \label{eq 5.1}%
\end{equation}
for all compact set $e$\ with $cap\left(  e\right)  =cap\left(
e;H^{r}\right)  $. Suppose that, for all weights $w\in A_{1}$,
\begin{equation}
\int\limits_{\mathbb{R}^{d}}\left\vert g(x)\right\vert ^{2}wdx\leq
K\int\limits_{\mathbb{R}^{d}}\left\vert h(x)\right\vert ^{2}wdx \label{eq 5.2}%
\end{equation}
with a constant $K$ depending only on $d$ and the constant $A$ in the
Muckenhoupt condition. Then%
\[
\int\limits_{e}\left\vert g(x)\right\vert ^{2}dx\leq Ccap(e)
\]
for all compact sets $e$ with $C=C\left(  d,r,K\right)  $ .
\end{theorem}

To show this theorem, we need some facts from the equilibrium potential of a
compact set $e$ of positive capacity [AH]. The equilibrium potential of a
measure $\mu\in M^{+}$ is defined by%
\[
P=P_{e}=J_{r}\left(  J_{r}\mu\right)  .
\]

\begin{lemma}
[{[AH]}]\label{lemme2.2}For any compact set $e\subset$ $\mathbb{R}^{d}$, there
exists a measure $\mu=\mu_{e}$ such that

\begin{description}
\item[(i)] supp $\mu\subset e;$

\item[(ii)] $\mu(e)=cap\left(  e,H^{r}\right)  ;$

\item[(iii)] $\left\Vert J_{r}\mu\right\Vert _{L^{2}}^{2}=cap\left(
e,H^{r}\right)  ;$

\item[(iv)] $P_{e}(x)\geq1$ quasi-everywhere on $\mathbb{R}^{d};$

\item[(v)] $P_{e}(x)\leq K=K(d,r)$ on $\mathbb{R}^{d};$

\item[(vi)] $cap\left\{  P_{e}\geq t\right\}  \leq At^{-1}cap\left(
e,H^{r}\right)  $ for all $t>0$ and the constant is independent of $e$.
\end{description}
\end{lemma}

The measure $\mu_{e}$ associated with $e$ is called the capacitary
(equilibrium) measure of $e$. We will also need the asymptotics (Voir [AH])%
\begin{align}
G_{\alpha}(x)  &  \simeq\left\vert x\right\vert ^{\alpha-d}\text{ \ , \ if
}d\geq3\text{, \ \ }\left\vert x\right\vert \rightarrow0;\label{eq1.5}\\
G_{\alpha}(x)  &  \simeq\left\vert x\right\vert ^{\frac{\alpha-d}{2}%
}e^{-\left\vert x\right\vert }\text{ \ , \ if }d\geq2\text{, \ \ }\left\vert
x\right\vert \rightarrow+\infty;\nonumber
\end{align}
Sometimes, it will be more convenient to use a modified kernel
\[
\overset{\sim}{G_{r}}(x)=\max\left(  G_{r}(x),1\right)
\]
which does not have the exponential decay at $\infty$. Obviously, both $G_{r}$
and $\overset{\sim}{G_{r}}$ are positive nonincreasing radial kernels.
Moreover, $\overset{\sim}{G_{r}}$ has the doubling property :%
\[
\overset{\sim}{G_{r}}(2s)\leq\overset{\sim}{G_{r}}(s)\leq c(d)\overset{\sim
}{G_{r}}(2s)
\]
The corresponding modified potential is defined by%
\[
\overset{\sim}{P}(x)=\overset{\sim}{G_{r}}\ast\mu(x)
\]
The rest of the proof of theorem \ref{thMV} is based on the following
proposition :

\begin{proposition}
\label{poids}Let $d\geq2$ and let $0<\delta<\frac{d}{d-2r}$.\ Then
$\overset{\sim}{P}$ lies in the Muckenhoupt class $A_{1}$on $\mathbb{R}^{d}$,
i.e.,
\[
M\left(  \overset{\sim}{P}^{\delta}(x)\right)  \leq C(\delta,d)\overset{\sim
}{P}^{\delta}(x),\ dx\text{ p.p}%
\]
where $M$ denotes the Hardy-Littlewood maximal operator on $\mathbb{R}^{d}$,
and the corresponding $A_{1}-$bound $C(\delta,d)$ depends only on $d$ and
$\delta.$
\end{proposition}

\begin{proof}
Let $k:\mathbb{R}_{+}\rightarrow\mathbb{R}_{+}$ be a nonincreasing function
which satisfies the doubling condition :%
\[
k(2s)\leq ck(s)\text{, }s>0
\]
It is easy to see that the radial weightl $k(\left\vert x\right\vert )\in
A_{1}$ if and only if%
\begin{equation}
\int\limits_{0}^{R}k^{\delta}(t)t^{d-1}dt\leq cR^{d}k(R)\text{, }R>0
\label{eq5.6}%
\end{equation}
Moreover, the $A_{1}-$bound \ of $k$ is bounded by a constant which depends
only on $C$ in the preceding estimate and the doubling constant $c$ (see
[St2]). It follows from that
\[
\overset{\sim}{G}_{r}(s)\simeq\left\vert s\right\vert ^{r-d}\ \ \text{if
\ }d\geq3\text{ }\ \text{for }\ 0<s<1
\]
and
\[
\overset{\sim}{G}_{r}(s)\simeq1\ \text{for }s\geq1
\]
Hence, $k\left(  \left\vert s\right\vert \right)  =\overset{\sim}{G}%
_{r}^{\delta}(s)$ is a radial nonincreasing kernel with the doubling property.
By Jensen' inequality, we have
\[
\overset{\sim}{G}_{r}^{\delta_{1}}\in A_{1}\text{ implies }\overset{\sim}%
{G}_{r}^{\delta_{2}}\in A_{1}\text{ \ if \ }\delta_{1}\geq\delta_{2}%
\]
Clearly (\ref{eq5.6}) holds if and only if $0<\delta<\frac{d}{d-2r}$. Hence,
without loss of generality, we assume $1\leq\delta<\frac{d}{d-2r}$. Then by
Minkowski's inequality and the $A_{1}-$estimate for $\overset{\sim}{G}%
_{r}^{\delta}$ established above, it follows%
\begin{align*}
M\left(  \overset{\sim}{P}^{\delta}(x)\right)   &  \leq M\left(  \left(
\overset{\sim}{G}_{r}^{\delta}\right)  ^{\frac{1}{\delta}}\ast\mu(x)\right)
^{\delta}\\
&  \leq C(\delta,d)\left(  \overset{\sim}{G}_{r}\ast\mu\right)  ^{\delta}(x)\\
&  =C(\delta,d)\overset{\sim}{P}^{\delta}(x).\text{ }%
\end{align*}

\end{proof}

We are now in a position to prove theorem \ref{thMV}.

\begin{proof}
Suppose $\upsilon_{e}$ is the capacitary measure of $e\subset\mathbb{R}^{d}$
and let $\varphi=P$ is its potential. Then, by lemma \ref{lemme2.2}, we have

\begin{description}
\item[(i)] $\varphi(x)\geq1$ quasi-everywhere on $e$ ;

\item[(ii)] $\varphi(x)\leq B=B\left(  d,r\right)  $ for all $x\in
\mathbb{R}^{d}$ ;

\item[(iii)] $cap\left\{  \varphi\geq t\right\}  \leq Ct^{-1}cap\left(
e\right)  $ for all $t>0$ with the constant $C$ is independent of $e$.

\item Now, it follows from a proposition \ref{poids} that $\varphi^{\delta}\in
A_{1}.$ Hence, by (\ref{eq 5.2}),
\[
\int\limits_{\mathbb{R}^{d}}\left\vert g(x)\right\vert ^{2}\varphi^{\delta
}dx\leq K\int\limits_{\mathbb{R}^{d}}\left\vert h(x)\right\vert ^{2}%
\varphi^{\delta}dx
\]
Applying this together with $(i)$ and $(ii)$, we get%
\[
\int\limits_{e}\left\vert g(x)\right\vert ^{2}dx\leq\int\limits_{\mathbb{R}%
^{d}}\left\vert g(x)\right\vert ^{2}\varphi^{\delta}dx\leq C\int
\limits_{\mathbb{R}^{d}}\left\vert h(x)\right\vert ^{2}\varphi^{\delta
}dx=C\int\limits_{0}^{B}\int\limits_{\varphi\geq t}\left\vert h(x)\right\vert
^{2}dxt^{\delta-1}dt
\]
By $(\ref{eq 5.1})$ and $($iii$)$,
\[
\int\limits_{\varphi\geq t}\left\vert h(x)\right\vert ^{2}dx\leq Ccap\left\{
\varphi\geq t\right\}  \leq\frac{C}{t}cap\left(  e\right)
\]
Hence,%
\[
\int\limits_{e}\left\vert g(x)\right\vert ^{2}dx\leq C\int\limits_{0}%
^{B}t^{-1}cap\left(  e\right)  t^{\delta-1}dt=Ccap\left(  e\right)
\int\limits_{0}^{B}t^{\delta-2}dt
\]
Clearly, for all $0\leq r<\frac{d}{2}$, we can choose $\delta>1$ so that
$0<\delta<\frac{d}{d-2r}$. Then
\[
\int\limits_{0}^{B}t^{\delta-2}dt=\frac{B^{\delta-1}}{\delta-1}<\infty
\]
which concludes
\[
\int\limits_{e}\left\vert g(x)\right\vert ^{2}dx\leq Ccap\left(  e\right)  .
\]

\end{description}
\end{proof}

\end{document}